\documentclass[a4paper,12pt,reqno]{article}
\usepackage[utf8]{inputenc}

\usepackage{leftidx}
\usepackage{amsfonts}
\usepackage{amsmath}
\usepackage{amssymb}
\usepackage{amsthm}
\usepackage{enumitem}
\usepackage[T1]{fontenc}
\usepackage{cite}
\usepackage[top=20mm,bottom=20mm,left=20mm,right=20mm]{geometry}
\usepackage[pdfstartview=FitV,bookmarks=false]{hyperref}
\usepackage{cite}
\usepackage{tikz}
\usepackage{float}

\usepackage[french]{babel}

\newmuskip\pFqmuskip

\newcommand*\pFq[6][8]{%
  \begingroup % only local assignments
  \pFqmuskip=#1mu\relax
  \mathchardef\normalcomma=\mathcode`,
  \mathcode`\,=\string"8000
  \begingroup\lccode`\~=`\,
  \lowercase{\endgroup\let~}\pFqcomma
  {}_{#2}F_{#3}{\left[\genfrac..{0pt}{}{#4}{#5};#6\right]}%
  \endgroup
}
\newcommand{\pFqcomma}{{\normalcomma}\mskip\pFqmuskip}

\allowdisplaybreaks

\font\Bbb=msbm10

\def\S{\hbox{\Bbb S}}
\def\D{\hbox{\Bbb D}}

\newcommand{\C}{\ensuremath{\mathbb{C}}}

\newcommand{\mathP}{\ensuremath{\mathbb{P}}}
\newcommand{\Q}{\ensuremath{\mathbb{Q}}}

\newcommand{\Z}{\ensuremath{\mathbb{Z}}}
\def\img{{\bf i }}

\def\b0{{\bf 0}}

\def\bj{{\bf j}}
\def\bn{{\bf n}}

\def\b1{{\bf 1}}

\def\bGamma{\bf \Gamma}

\def\bh{{\bf h}}

\def\bj{{\bf j}}
\def\bk{{\bf k}}

\def\bn{{\bf n}}

\def\b1{{\bf 1}}

\def\Q{\hbox{\Bbb Q}}
\def\S{\hbox{\Bbb S}}
\def\C{\hbox{\Bbb C}}
\def\D{\hbox{\Bbb D}}
\def\J{\hbox{\Bbb J}}

\def\Z{\hbox{\Bbb Z}}

\numberwithin{equation}{section}

\newtheorem{theorem}{\bf Th\'eor\`eme}[section]
\newtheorem{lemma}[theorem]{\bf  Lemme}
\newtheorem{remark}[theorem]{\bf Remarque}

\newtheorem{corollary}[theorem]{\bf Corollaire}
\newtheorem{definition}[theorem]{\bf D\'efinition}
\newtheorem{conjecture}[theorem]{\bf Conjecture}

\date{}
\title {
Note sur le lieu discriminant d'une  hypersurface de bi-degr\'e $(m, n).$   } 
\author{ Susumu TANAB\'E} 
\begin{document}
\maketitle
%\noindent
%

\begin{flushright}
En hommage \`a Valentin Po\'enaru
\end{flushright}

\selectlanguage{french}
\begin{abstract}
Nous pr\'esentons une  m\'ethode topologique pour l'\'etude du lieu discriminant d'une vari\'et\'e alg\'ebrique d\'efinie dans un produit d'espaces projectifs. Notre approche se distingue par l'utilisation efficace du groupo\"ide dans la description de la monodromie.
Comme exemple, nous traitons ici le lieu discriminant d'une hypersurface de bi-degr\'e $(m, n).$
\end{abstract}

\section{Groupo\"ide \`a la recherche  de la monodromie \label{sec:groupoide}}

Dans ce chapitre, nous pr\'esentons une m\'ethode topologique pour l'\'etude du lieu discriminant d'une vari\'et\'e alg\'ebrique.
Vu la propri\'et\'e remarquable du lieu discriminant d'une vari\'et\'e d\'efinie dans un produit d'espaces projectifs, mentionn\'ee dans la remarque \ref{covering} ci-dessous, nous esp\'erons que notre approche pourra contribuer \`a l'\'eclaircissement de certaines propri\'et\'es g\'eom\'etriques et analytiques cach\'ees au fond du ph\'enom\`ene de monodromie. En tant qu'objectif d'importance actuelle et pourtant r\'ealisable, nous pouvons mentionner l'application \`a la sym\'etrie miroir homologique \cite[corollary 3.5, theorem 3.6]{Horja}, \cite{Tan17}. 
Dans le but de simplifier l'exposition, on se restreint ici au cas d'une hypersurface de bi-degr\'e $(m, n)$ \cite{TanLonne}, surtout au cas $n=2.$

Nous esp\'erons que les exemples que nous pr\'esentons dans ce chapitre et dans les autres articles en pr\'eparation (\cite{Tan24}, \cite{TanKocar}, \cite{TanLonne}), illustrent la th\`ese promue par A. Grothendieck \cite{Grothendieck}
sur l'avantage des groupo\"ides dans l'\'etude du groupe fondamental du compl\'ement du lieu discriminant.

{\it M\^eme pour $g=0$ (donc quand les groupes de Teichm\"uller correspondants sont des groupes de tresses ''bien connus''  ),
les g\'en\'erateurs et relations connus \`a ce jour dont j'ai eu connaissance, me semblent inutilisables tels quels, car ils ne pr\'esentent pas les caract\`eres
d'invariance et de sym\'etrie indispensables pour que l'action de $\bGamma = Gal (\bar Q/Q)$ soit directement lisible sur cette pr\'esentation.  
Ceci est li\'e notamment au fait que les gens s'obstinent encore, en calculant avec des groupes fondamentaux, \`a fixer un seul point base, plut\^ot que d'en choisir
astucieusement tout un paquet qui soit invariant par les sym\'etries de la situation. Dans certaines situations (comme des th\'eor\`emes de descente \`a la Van Kampen pour groupes fondmentaux )
il est bien plus \'el\'egant, voir indispensable pour y comprendre quelque chose, de travailler avec des groupo\"ides fondamentaux par rapport \`a un paquet de points base convenable,
et il en est certainement ainsi pour la tour de Teichm\"uller.
}

\footnoterule

\footnotesize{AMS Subject Classification:  20F36  (primary), 14D05, 32S40
(secondary).
Key words and phrases: discriminantal loci, braid monodromy,  groupoid.}
\footnotesize{Partiellement soutenu par  Max Planck Institut f\"ur Mathematik, Research fund of the Galatasaray University project number 1123 "Analysis and topology of algebraic functions and period integrals."}
\normalsize

\newpage
On trouve la description de la m\^eme th\`ese exprim\'ee d'une mani\`ere l\'eg\`erement modifi\'ee dans \cite[3. Application of fundamental groupoid]{Brown}.  

Notamment, dans \cite{Tan24}  nous proposons la recherche de la monodromie globale des int\'egrales de p\'eriodes associ\'ees 
\`a une hypersurface (ou une intersection compl\`ete)  affine (vari\'et\'e Calabi--Yau)  en se servant de leur prolongement analytique le long de chemins qui \'evitent 
les lieux discriminants. On peut, donc, interpr\'eter notre essai comme une sorte de continuation dans l'esprit de la th\`ese ci-dessus 
o\`u le centre d'attention se d\'eplace vers la monodromie de la vari\'et\'e alg\'ebrique affine 
du groupe de Galois absolu  $\bGamma = Gal (\bar Q/Q).$

Notre programme propos\'e dans \cite{Tan24} a \'et\'e r\'ealis\'e pour la courbe de bi-degr\'e $(2, 2)$ d\'efinie dans $\mathP^1 \times \mathP^1$  \cite{Tan17}.

Une analyse de la connexion \`a singularit\'es r\'eguli\`eres et irr\'eguli\`res \cite{PaulRamis} m\'erite une mention sp\'eciale.  
Cet exemple  illustre \'eloquemment
la mise en valeur de la th\`ese ci-dessus dans le cadre de la recherche globale des \'equations diff\'erentielles.

Il faut signaler que le contenu du chapitre actuel sert d'\'etape pr\'eliminaire pour la r\'ealisation d'un projet pr\'epar\'e avec le concours de Michael L\"onne \cite{TanLonne}.  L'auteur remercie Mutlu Ko\c{c}ar de l'assistance \`a la pr\'eparation des figures et Athanase Papadopoulos de la lecture soigneuse du texte.  

\section{Lieu discriminant d'une hypersurface de bi-degr\'e $(m, n)$ \label{sec:bidegree}}

%\subsection{Uniformisation de Horn-Kapranov pour le lieu discriminant }5

Cette section est consacr\'ee \`a l'analyse de l'espace des modules
d'une hypersurface $Y$ qui est la vari\'et\'e de
sym\'etrie miroir d'une vari\'et\'e Calabi--Yau de bi-degr\'e $(m, n)$ g\'en\'erique 
dans $\mathP^{m-1} \times \mathP^{n-1},$
au sens de V.V. Batyrev  
 \cite{Baty3}, \cite[4.1]{CK}.  D'apr\`es Batyrev, 
le poly\`edre de  Newton $\Delta(Q)$ du polyn\^ome de Laurent qui d\'efinit la partie affine de $Y$
$$ Q(z,w)= \frac{P(z,w)}{z_1...z_{m-1} w_1...w_{n-1}} $$
est  le dual polaire (polar dual) de l'enveloppe convexe d'un produit de deux simplexes canoniques de dimensions $m-1$ et $n-1$.

En vertu de l'action des tores alg\'ebriques sur les variables $(z_1,...,z_{m-1}, w_1, ..., w_{n-1}),$
le polyn\^ome de Laurent $Q(z,w)$ \`a coefficients g\'en\'eraux se transforme en 
\begin{flushleft}
$  F_{m,n}(z,w) =$
\end{flushleft}
\begin{equation}
1 + z_1 +z_2 + \cdots + z_{m-1}+\frac{x}{z_1\cdots z_{m-1} } + w_1 + \cdots +w_{n-1} +\frac{y}{w_1\cdots w_{n-1}}
\label{Fnm}
\end{equation}
qui d\'epend de deux param\`etres $(x,y)$  de d\'eformation. 

 Le lieu discriminant $S$ de la vari\'et\'e affine
\begin{equation}
X^{(m,n)}_{x,y}= \{(z,w) \in (\C^{\ast})^{m+n-2} ;  F_{m,n}(z,w)=0 \}
\label{Xmn}
\end{equation}
se pr\'esente  comme la r\'eunion d'une courbe 
\begin{equation}
C=\{(x,y); D(x,y)=0\}
\label{Ccourbe}
\end{equation}
et des deux axes de coordonn\'ees,
\begin{equation}S=\{(x,y) \in \C^2 ; x y D(x,y)=0 \}.
\label{Sdiscriminant}
\end{equation}
D'apr\`es la d\'efinition du lieu discriminant, la vari\'et\'e affine $X^{(m,n)}_{x,y}$ \eqref{Xmn} est singuli\`ere si et seulement si $(x,y) \in S.$

La courbe  $C$ \eqref{Ccourbe} admet l'uniformisation de Horn--Kapranov \cite{Tan07}
comme suit
\begin{equation} (x,y) = ((\frac{s}{ms +n})^m,  (\frac{1}{ms + n})^n  ),\; s \in \C
\label{uniformisation}
\end{equation}
de telle sorte que
$$ D\bigl((\frac{s}{ms +n})^m, (\frac{1}{ms + n})^n \bigr) =0 ,\; \forall s \in \C. $$

La d\'erivation de l'uniformisation \eqref{uniformisation} est bas\'ee sur le calcul du lieu singulier des int\'egrales de p\'eriodes
de la vari\'et\'e \eqref{Xmn}. Il est connu \cite{Berglund} que l'int\'egrale de p\'eriodes
$$ f(x,y)= \int_\gamma \frac{dz \wedge dw}{d F_{m,n}(z,w)},$$
pour le cycle quelconque $ \gamma \in H_{m+n-3} ( X^{(m,n)}_{x,y})$ satisfait au syst\`eme des \'equations hyperg\'eom\'etriques du type de Horn ci-dessous,
\begin{equation}
\begin{array}{l}
\left(\theta_{x}^m  - x (m\theta_x + n\theta_y+1)\cdots (m\theta_x + n\theta_y+m)\right)f(x,y)=0, \\
\left(\theta_{y}^n  - y (m\theta_x + n\theta_y+1)\cdots (m\theta_x + n\theta_y+n)\right)f(x,y)=0,
\end{array}
\label{simplicialhorn}
\end{equation}
pour $m \geq n \geq 2,$ 
$$ \theta_x = x \frac{\partial }{\partial x}, \theta_y = y \frac{\partial }{\partial y}.$$ 
La dimension de l'espace vectoriel complexe des solutions \eqref{simplicialhorn}, nomm\'ee rang holonomique,  
est \'egale \`a $mn$ $= |\chi(X^{(m,n)}_{x,y})|.$
Une base explicite de l'espace  
des solutions de \eqref{simplicialhorn} peut \^etre donn\'ee par les int\'egrales de Mellin--Barnes,
\begin{equation}
 f^{k}_{i,j}(x,y)= \int_{\Delta }\frac{\Gamma(s)^{m-i}\Gamma(t)^{n-j}\Gamma(1-(ms+nt))}{\Gamma(1-s)^i\Gamma(1-t)^j}x^{-s}y^{-t} e^{-\pi {\sqrt -1}(is+jt)}ds \wedge dt,
\label{ex(1,0)(0,1)(2,2)}
%\label{simplicialhorn}
\end{equation}
o\`u $( i,j) \in [0; m-1] \times [0; n-1]$ et $\Delta$ est  un semi-groupe form\'e des p\^oles de l'int\'egrand de  \eqref{ex(1,0)(0,1)(2,2)},
$$\Delta= \{ (s,t) \in \C^2 ; s \in \Z_{\leq 0 }, t \in \Z_{\leq 0} \}.$$
Ainsi  le lieu de ramification des solutions des \'equations diff\'erentielles 
\eqref{simplicialhorn} admet la repr\'esentation \eqref{uniformisation}, \cite{Tan07}.

Le groupe de monodromie globale des int\'egrales 
\eqref{ex(1,0)(0,1)(2,2)} s'obtient comme repr\'esentation de $\pi_1(\C^2 \setminus S).$ Le sujet de notre chapitre
fait partie d'un programme de recherche sur ce groupe de monodromie globale.

A partir d'ici,  on  utilise la notation
$\omega_m =e^{2\pi \sqrt -1/m},  \omega_n^j = e^{2\pi j \sqrt -1/n}$
ainsi que $e(\lambda)= e^{2\pi \sqrt -1  \lambda},$
c.\`a.d. $\omega_m^i =e(i/m), \omega_n^j =e(j/n).$ 
Nous utilisons aussi la notation $i \in [r_1; r_2] \Leftrightarrow i \in \{r_1, \cdots,  r_2\} $ pour deux entiers $r_1 \leq r_2.$
 
Afin d'\'etudier  le groupe fondamental $\pi_1 (\C^2 \setminus S),$ on commence par le lemme suivant.

\begin {lemma}
a) La courbe $C=\{(x,y); D(x,y)=0\}$ est irr\'educible.
b)La courbe $C$ est nodale c.\`a.d.  qu'elle poss\`ede seulement des points singuliers o\`u les deux branches se croisent
transversalement. c) $\pi_1(\C^2 \setminus C) \cong \Z.$
\label{irreduciblecurve}
\end {lemma}

\begin{proof}[D\'emonstration]
a) Gr\^ace  \`a \eqref{uniformisation} la branche $(x, y_j(x)), j=1, \cdots, m$ de la courbe $C$ 
admet le d\'eveloppement de Puiseux suivant, 
$$ y_j(x)= (\frac{1- m \omega_m^j x^{1/m}}{n})^n.$$
Puisque chaque branche  $y_j(x)$ appartient \`a la m\^eme classe de conjugaison, et puisque l'irr\'eductibilit\'e d'une courbe
est \'equivalente \`a l'unicit\'e de la classe de conjugaison des branches, on a l'\'enonc\'e. 

b) Il suffit de d\'emontrer que la courbe ayant la pr\'esentation param\'etrique
$C'=\{(x,y) ; x=(1-t)^n, y=t^m \}$
ne poss\`ede que des points singuliers doubles. 
Consid\'erons les polyn\^omes $P_{k,\ell}(t) = (1-\omega_n^k t)^m -(1-\omega_n^\ell t)^m$ pour $0 \leq  k \not = \ell \leq n-1.$

Nous avons besoin de d\'emontrer que (1)  $P_{0,\ell}(t)=0$, $ \ell \not= 0$ ne poss\`ede aucune racine multiple.  (2) Le syst\`eme 
d'\'equations $P_{0,\ell}(t)= P_{k,\ell}(t)=0 $ n'admet aucune solution non nulle pour
 $0 <k \not = \ell \leq m-1.$

Pour que l'\'equation $P_{0,\ell}(t)=0$ admette une racine double, l'existence  d'une  solution du syst\`eme suivant est n\'ecessaire: $ P_{0,\ell}(t)=  P'_{0,\ell}(t)=0.$ Si $t_0$ est une solution du syst\`eme, alors elle doit satisfaire $\omega_n^\ell P_{0,\ell }(t_0) +(1- \omega_n^\ell t_0)P'_{0,\ell}(t_0)/m = 0 ,$ c. \`a. d.  $(1-t_0)^{m-1}(\omega_n^\ell-1)=0.$ Puisque $P_{0,\ell}(1) = -(1-\omega_n^\ell)^m \not =0,$ nous avons (1).

Les racines non nulles de  $P_{0,\ell}(t)=0$ sont de la forme $\{\frac{1- \omega_m^i}{1- \omega_m^i\omega_n^\ell}\},$
$1 \leq i \leq m-1,  (\frac{\ell}n + \frac{i}{m}) \not \equiv 0 \;  mod \;1,$ pourtant celles de  $P_{k,\ell}(t)=0 $ sont $\{\frac{1- \omega_m^j}{\omega_n^k- \omega_m^j\omega_n^\ell}\},$ $1\leq  j \leq m-1,$    $(\frac{-k+\ell}{n} + \frac{j}{m}) \not \equiv 0 \; mod \;1.$  On rappelle ici que $t=0$ correspond \`a un point  $(m-1)-$flexe de la courbe $C'.$ L'argument $Arg \left(\frac{1- \omega_m^j}{\omega_n^k- \omega_m^j\omega_n^\ell}/\frac{1- \omega_m^i}{1- \omega_m^i\omega_n^\ell}\right) = -\frac{k}{n}\pi $
ne devient jamais  $0$ ou $2 \pi$ pour $1 \leq k \leq n-1.$  Ceci montre (2).

c) D'apr\`es ~\cite{Orevkov88}, la condition suivante (appel\'ee "negativity condition at infinity": abbr\'eviation NC) 
est suffisante pour qu'une courbe nodale  $C$ poss\`ede un groupe ab\'elien $\pi_1(\C^2 \setminus C)$.

 \vspace{1pc}

({\bf NC})
 La restriction de la fonction alg\'ebrique $y(x)$ satisfaisant $D(x,y)=0$ au voisinage de l'infini (l'axe r\'eel exclu) consiste en $n$ branches analytiques $y_j(x)$ de la  forme $g_j(\tau(x))$ o\`u $g_j(\tau)$ est une fonction m\'eromorphe et  
univalente au voisinage du  point  $\tau=0$ pour $\tau(x):$ une branche univalente de la fonction $x^{-1/n}.$
On a $\lim_{\tau \rightarrow 0} (g_j(\tau)- g_k(\tau)) \not =0$ pour $k\not =j$  c.\`a.d.  qu'elle y poss\`ede soit un p\^ole, soit une limite finie.

\vspace{1pc}

Il s'ensuit de \eqref{uniformisation} et d'un changement lin\'eaire de variables que la fonction 
$g_j(\tau) = (\frac{1}{n})^n(1-\frac{m \omega_m^j}{ \tau} )^n$  satisfait {\bf NC}. 

Par le th\'eor\`eme de Hurewicz on a $\pi_1(\C^2 \setminus C)/Z = H_1(\C^2 \setminus C),$
mais la discussion ci-dessus montre que le commutateur $Z$ s'annule. Il est bien connu que $H_1(\C^2 \setminus C) \cong \Z$ pour la courbe irr\'eductible $C.$ 
\end{proof}

Ainsi $\pi_1(\C^2 \setminus S)$ est engendr\'e par trois g\'en\'erateurs $\gamma_1, \gamma_2, \gamma_C$
o\`u $\gamma_1$ (resp. $ \gamma_2,$  $\gamma_C$)  correspond au lacet qui fait un tour autour de la courbe $ x =0$ (resp. $ y=0,$   $ C$).

L'uniformisation \eqref{uniformisation} indique que le rev\^etement  de degr\'e $m n$ du lieu discriminant \eqref{Sdiscriminant} est homotope \`a l'arrangement complexe de $(mn+2)$ droites.
\begin{equation}
{\mathcal A}_{m,n} :=\{(s,t) \in \C^2; s t \prod_{0 \leq i \leq m-1, 0 \leq j \leq n-1}  ( \omega_m^i s+\omega_n^j t +1)=0 \}.
\label{Amn}
\end{equation}

La suite exacte suivante est une cons\'equence du th\'eor\`eme de Reidemeister--Schreier\cite{Zieschang}
\begin{equation}
 1 \rightarrow \pi_1(\C^2 \setminus  {\mathcal A}_{m,n}) \rightarrow \pi_1(\C^2 \setminus  S) \rightarrow^{h} (\Z/m\Z) \times
(\Z/m)  \rightarrow 1. 
\label{exactsequence3} 
\end{equation}
avec $ h(\gamma_C) = h(\gamma_1^m) =  h(\gamma_2^n)=1.$

\begin{remark}
{\rm
En effet, c'est un exemple de ph\'enom\`enes plus g\'en\'eraux
qu'on observe sur le lieu discriminant d'une hypersurface Calabi--Yau $X_\bn$ de multi-degr\'e $\bn =(n_i)_{i=1}^K$ d\'efinie dans le produit
d'espaces projectifs $\prod_{i=1}^K \mathP^{n_i-1}.$  
Le calcul de \cite[section 8]{Stienstra1} indique que le rev\^etement  de degr\'e $\prod_{i=1}^K {n_i}$ du lieu discriminant d'une telle hypersurface
admet la pr\'esentation comme arrangement d'hyperlans
\begin{equation}
{\mathcal A}_{\bn} :=\{(t_1, \cdots, t_K)  \in \C^K; \prod_{i=1}^K t_i  \cdot \prod_{(j_1, \cdots, j_K)\in \prod_{i=1}^K (\Z/n_i\Z)}
 \bigl( \sum_{\ell=1}^K \omega^{j_\ell}_{n_\ell} t_\ell+1 \bigr)=0 \}
\label{Abn}
\end{equation}
qui consiste de $ K+ \prod_{i=1}^K n_i$ hyperplans. La d\'emonstration se fait \`a l'aide des formules donn\'ees dans
\cite[theorem 2.6]{Tan07}.

Pour le cas o\`u $n_i = 2, \forall i \in [1; K],$ T. Terasoma \cite{Terasoma} utilise le rev\^etement 
\eqref{Abn} afin d'obtenir une description compl\`ete du groupe fondamental de l'ensemble des hypersurfaces lisses
de multi-degr\'e $\bn =(2,\cdots, 2).$
 }
\label{covering}
\end{remark}

\section{Lieu discriminant au cas $n=2.$}

\subsection{Rotation des poin\c{c}ons et tresse annulaire}
 
Dans cette section, nous \'etudions le groupe fondamental $\pi_1(\C^2 \setminus  {\mathcal A}_{m,2})$
 au moyen de la monodromie de tresses pour $m \geq 2$.

Tout d'abord, on calcule les branches $s(t)$ de l'arrangement ${\mathcal A}_{m,2}$ \eqref{Amn} comme fonctions complexes lin\'eaires de $t.$ On met
\begin{equation}
T_{j,\ell} (t)=- \omega^j ( 1+(-1)^\ell t ),
\label{Tjl}
\end{equation}
pour $j \in [0; m-1], \ell \in [0;1]$ et $\omega = \omega_m =e(1/m).$ D\'esormais  nous consid\'erons toujours  les indices $(j, \ell)$ dans $\Z/m \times \Z/2.$
%Nous utilisons aussi la notation $i \in [r_1; r_2] \Leftrightarrow i \in \{r_1, \cdots,  r_2\} $ pour deux entiers $r_1 \leq r_2.$

Nous suivons le changement de la section de ${\mathcal A}_{m,2}$ en fonction de la valeur $t \in \C,$
\begin{equation}
\S(t) := \{0\} \bigcup_{[j,\ell ] \in  [0;m-1]\times [0;1] }  \{T_{j,\ell} (t) \}.
\label{sectionSt}
\end{equation}

Dor\'enavant nous nommons les \'el\'ements de $\S(t)$ {\bf poin\c{c}ons} (punctures), par contre,  la valeur fixe de $t$ sera nomm\'ee  {\bf point}.

On note 
\begin{equation}
\S^\ast(t) := \S(t) \setminus \{0\}.
\label{sectionStast}
\end{equation}

\begin{remark}
{\rm La discussion dans la section \ref{tressesinduites}  ci-dessus  montre que l'on peut rencontrer une difficult\'e lors de la g\'en\'eralisation aux cas $n \geq 3. $ En particulier,  dans les cas $n \geq 5$ intervient le ph\'enom\`ene d'entrelacement (entanglement) des branches $s(t)$
de l'arrangement ${\mathcal A}_{m,n}.$  Voir \cite[6.1 Parallel transport in the model family]{Loenne}.  Si cette difficult\'e n'exige pas un changement drastique de notre vision de la topologie de l'arrangement ${\mathcal A}_{m,n},$ nous pourrions affirmer un \'enonc\'e similaire \`a la conjecture \ref{conj} formul\'ee \`a la fin du chapitre.}
\end{remark}

\begin{lemma}\label{cooincidence}

Les branches \eqref{Tjl}  se rencontrent
aux points de co\"incidence  $\D \cup \{0\}$ avec
\begin{equation}
   \D : =\{   \pm 1\}\cup \{ \img \tan \frac{ j \pi}{m}  \}_{ j \in [0; m-1] \setminus \{ \frac m2\}},
\label{coincidenceD}
\end{equation} ayant les types suivants.
Dor\'enavant on utilise la notation 
\begin{equation}
  |\bar \bj| = j_2 - j_1,  \; |\bj| = j_1 + j_2,
\label{jell}
\end{equation}
pour $\bj = (j_1, j_2).$

(1) Les cas $|\bar \ell|$
 =0.  En $t = (-1)^{\ell+1},$ confluent les $m$ racines suivantes:
\begin{equation}
T_{0,\ell} (  (-1)^{\ell+1} ) = T_{1,\ell} (  (-1)^{\ell+1} ) = \cdots =  T_{m-1,\ell} (  (-1)^{\ell+1} ),
\label{tpm1}
\end{equation}
pour $ \ell \in [0;1].$

(2)  Les cas $|\bar \ell| $ =1.   En $t = \img \tan \frac{ |\bar \bj| \pi}{m}$ confluent les $2$ racines qui forment au total $m-$paires distinctes:
\begin{equation}
T_{j_1,0} ( \img \tan \frac{ |\bar \bj| \pi}{m}) = T_{j_2,1 } (  \img \tan \frac{ |\bar \bj| \pi}{m}) = - e(\frac{|\bj|}{2m})(\cos \frac{|\bar \bj|}{m}\pi )^{-1},
\label{tpm2}
\end{equation}
si $|\frac{|\bar \bj|}{m}| \not = \frac{1}{2}.$

%(b) En $t=0,$ confluent les $2$ racines qui forment au total $m-$paires distinctes:
%\begin{equation}
%T_{j,0} ( 0) = T_{j,1 } ( 0) = - \omega^j,
%\label{tpmj}
%\end{equation} $j \in [0; m-1].$
\end{lemma}

%{\bf Inversion du sujet. En ce point confluent les racines... Voir Grevisse, Bon usage, \S 379, No. 2. Dans la langue juridique, administrative, didactique, quand le sujet est une sorte de d\'efinition ou une \'enum\'eration. e.g. Sont ali\'en\'es un adolescent opprim\'e, une femme mal mari\'ee, un ouvrier astreint \`a un travail ingrat,...}

\begin{figure}[H]
\centering
\includegraphics[totalheight=14cm]{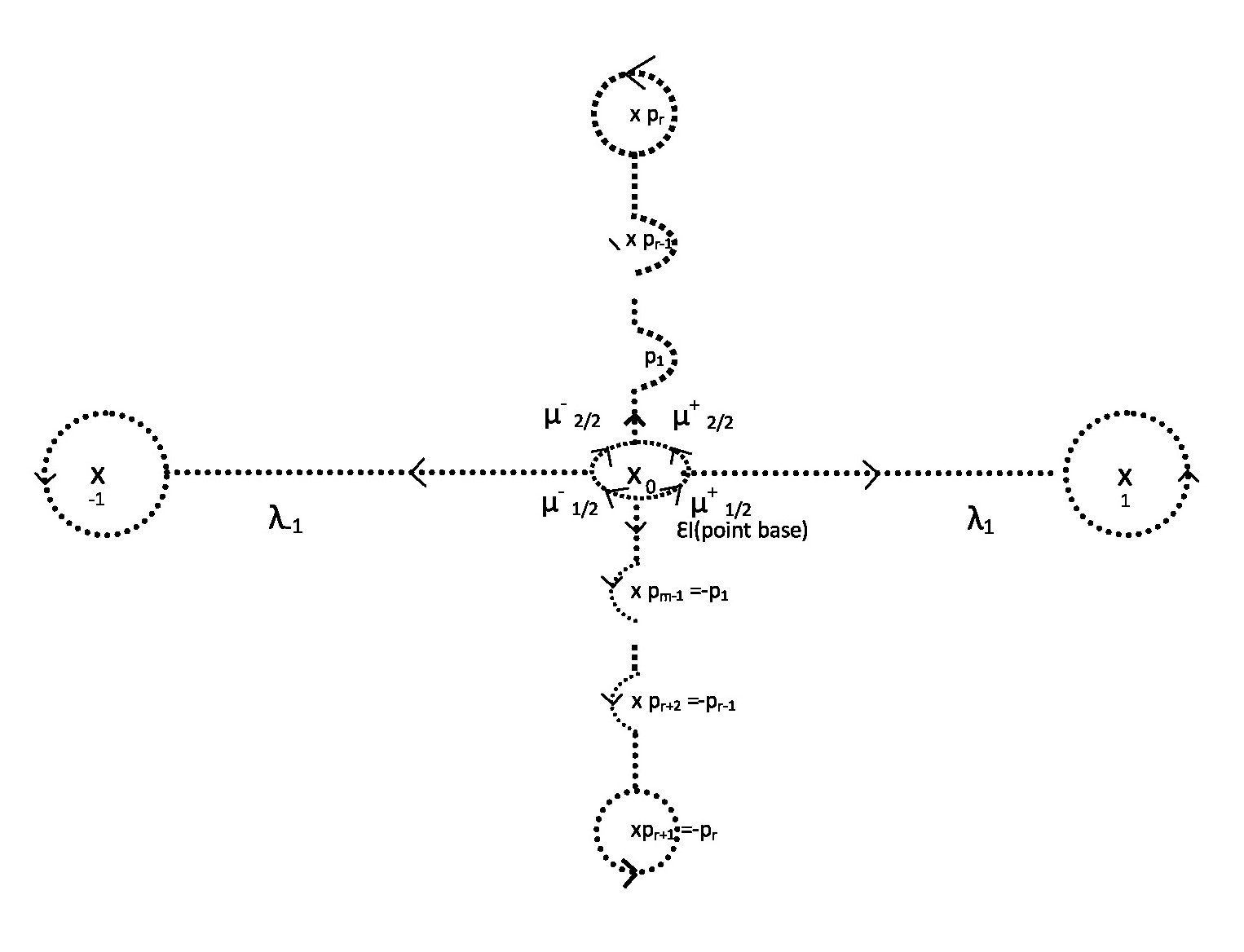}

\caption{ Figure de $\D$ avec des chemins issus du point base $t = - \epsilon \img.$}
\label{fig1}
\end{figure}

On appelle $p_k \in \D$ point de co\"incidence sur l'axe imaginaire pur,
\begin{equation}
p_k =  \img \tan \frac{ k \pi}{m},
\label{pk}
\end{equation}
pour $k \in [0;  m-1] \setminus \{\frac{m}{2}\}.$ Il s'ensuit que
\begin{equation}
 \tan \frac{ k \pi}{m}>0\;\; si\; k \in [1; r ],\; \tan \frac{ k \pi}{m}<0\; si \;  k \in [r+1; m-1 ],
\label{tanpos}
\end{equation} pour $r= \lfloor \frac{m-1}{2} \rfloor$, la partie enti\`ere de $\frac{m-1}{2}.$

Pour suivre le changement de la position des poin\c{c}ons 
\eqref{sectionSt},
nous fixons le point de d\'epart du mouvement de $t$ sur $\C^\ast \setminus \D.$
Par la sym\'etrie de la figure de $\D,$ nous choisissons $t = - \epsilon \img$
comme point base de d\'epart pour $\epsilon>0$ assez petit.

Au point base $t = - \epsilon \img,$ la fibre de r\'ef\'erence $\S ( - \epsilon \img)$ consiste en $2m$ poin\c{c}ons situ\'es sur un cercle
de centre $0.$

%{\bf Figure de $\S ( - \epsilon \img) $}

\begin{figure}[H]
\centering
\includegraphics[totalheight=12cm]{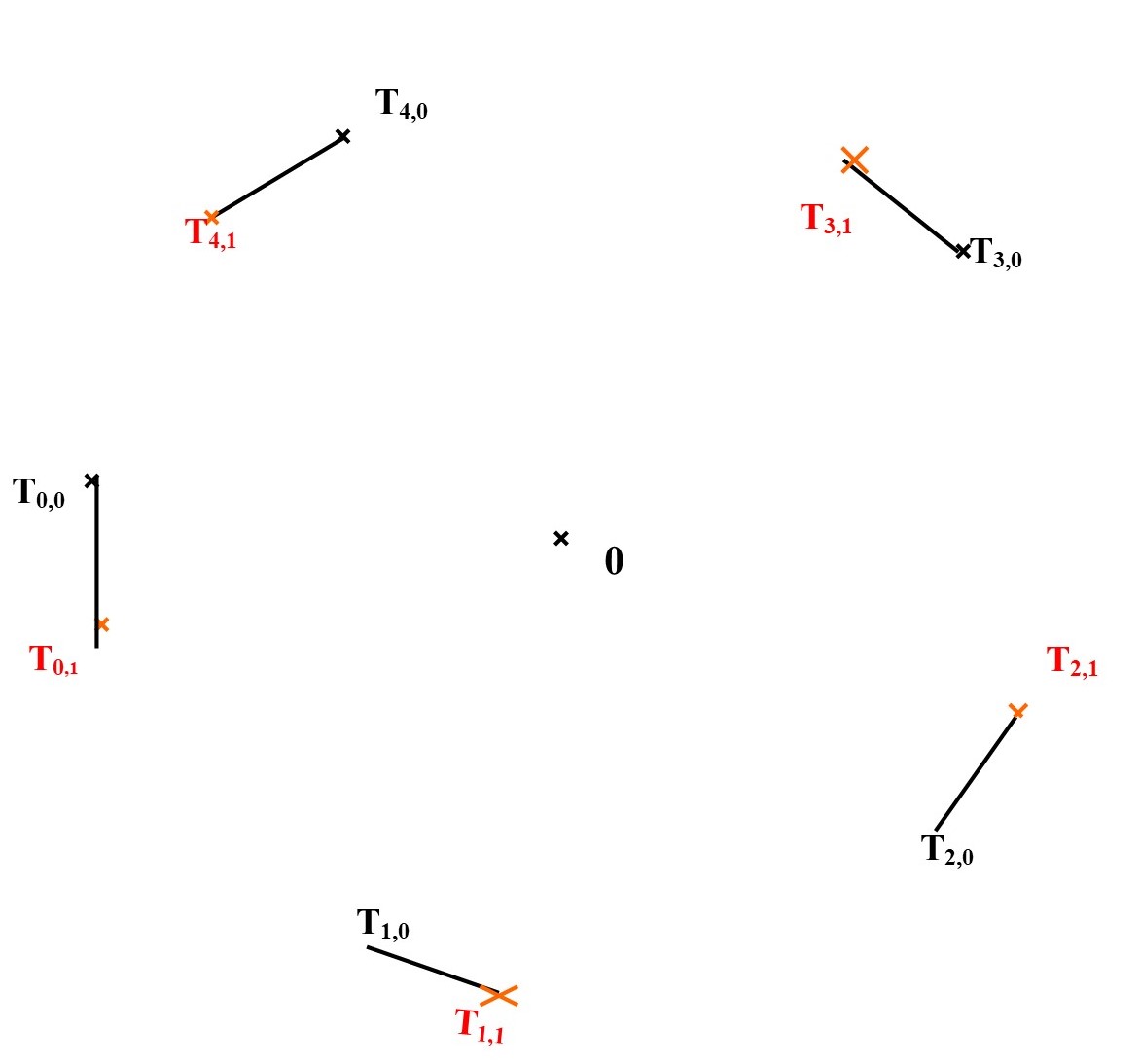}

\caption{ Figure de $\S ( - \epsilon \img) $ pour $m=5.$}
\label{fig2}
\end{figure}

Afin de d\'ecrire d'une mani\`ere efficace le mouvement de $\S(t)$ le long d'une courbe issue de $t = - \epsilon \img,$ 
nous nous rappelons la notion de {\bf rotation des poin\c{c}ons} (rotating punctures)\cite[Definition 1.4]{MM}.

\begin{definition}
Pour un ensemble d'indices $\J_\mu = \{(j_k, \ell_k)\}_{k=1}^\mu \subset [0, m-1] \times [0;1]$ et $t \in \C \setminus \D,$ on consid\`ere le polyg\^one $P_{\J_\mu}$
construit comme enveloppe convexe de $\mu$
 poin\c{c}ons $A_{\J_\mu}:=  \{ T_{j, \ell} (t)  \}_{(j, \ell) \in \J_\mu  }.$
% Pourtant on impose  la condition $ 0 \not \in  P_{\J_\mu},$ c.\`a.d.   $P_{\J_\mu}$ se trouve enti\`erement dans $\C^\ast.$
Supposons que les poin\c{c}ons $\{ T_{j_k, \ell_k} (t)  \}_{k=1}^\mu$ sont rang\'es selon le sens de l'orientation positive du bord  $\partial P_{\J_\mu}.$
Autrement dit, $\exists x \in( P_{\J_\mu})^{int} $ t.q. 
$$  Arg (T_{j_k, \ell_k} (t)  - x ) < Arg (T_{j_{k+1}, \ell_{k+1}} (t)  - x )   $$
$\forall (j_k, \ell_k) \in \J_\mu $
o\`u $\mu+1 \equiv 1.$

On d\'efinit la rotation $R_{A_{\J_\mu}}$ de l'ensemble  des racines  $A_{\J_\mu}$
comme le mouvement des racines ci-dessous,
\begin{equation}
R_{A_{\J_\mu}} (T_{j_k, \ell_k} (t)) = T_{j_{k+1}, \ell_{k+1}} (t).
\label{rotation}
\end{equation}
\label{defrotation}
\end{definition}

La puisssance de la rotation d\'ecoule imm\'ediatement de \eqref{rotation},
\begin{equation}
R_{A_{\J_\mu}}^p (T_{j_k, \ell_k} (t)) = T_{j_{\overline{k+p}}, \ell_{\overline{k+p}} }(t), 
\label{rotationpuissance}
\end{equation}
o\`u $\overline{k+p} \in [1; \mu]$ et $ k+p \equiv \overline{k+p} \; (mod \; \mu).$

\begin{figure}[H]
\centering
\includegraphics[totalheight=12cm]{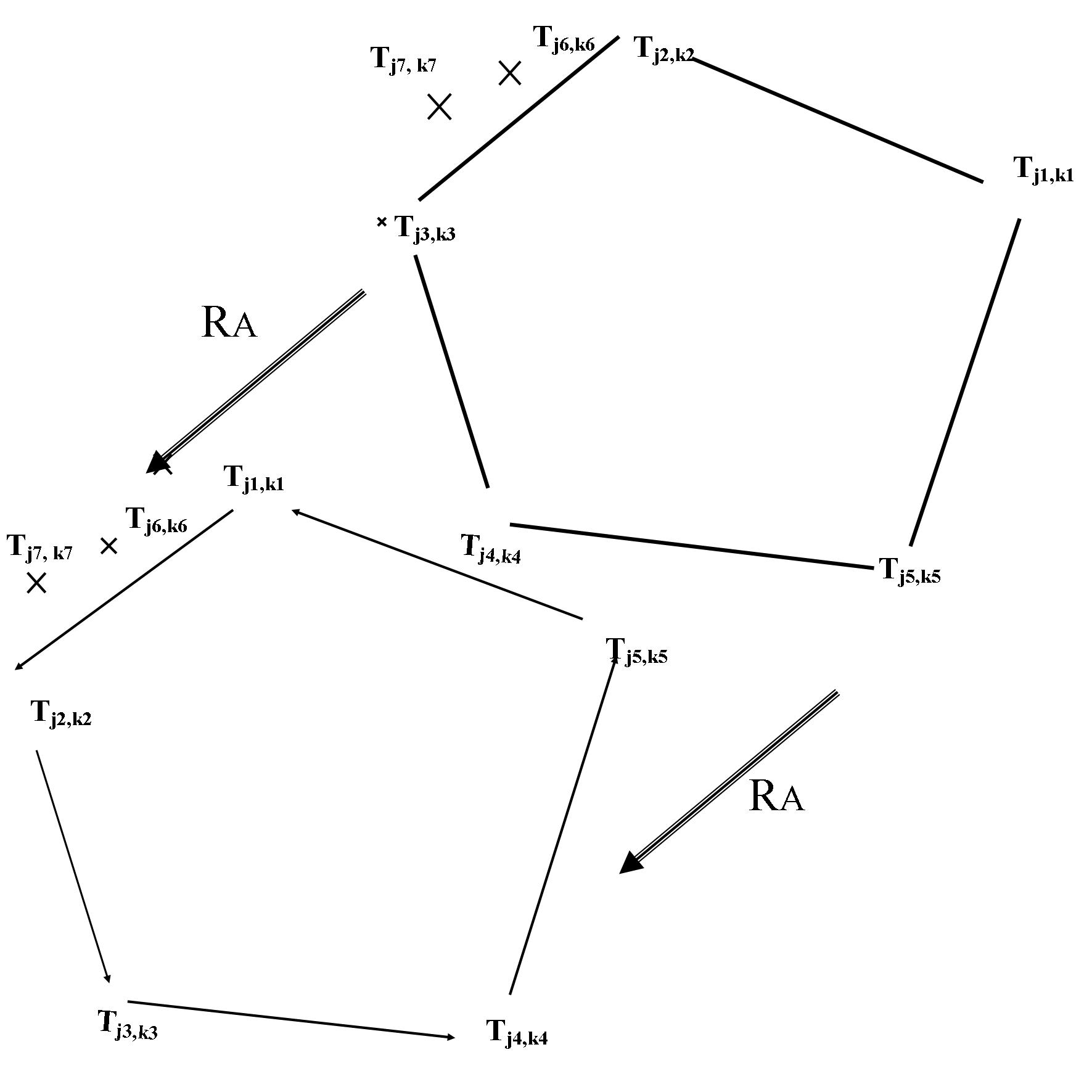}

\caption{Figure de $R_{A_{\J_\mu}}$ avec $\mu=5,$ $A=A_{\J_5} = \{T_{j_\ell,k_\ell}\}_{\ell=1}^5.$}
\label{figRAJ5}
\end{figure}

\begin{remark}
{\rm Si le mouvement de $t$ le long d'une courbe $ \lambda(u) \subset \C^\ast \setminus \D, u \in [0,1]$ t.q. $\lambda(0) = - \epsilon \img,$ provoque une rotation $R_{A_{\J}},$
alors nous obtenons une tresse \`a  $(2m+1)$ brins dans $\C \times [0,1]$ trac\'es par les racines $\S(\lambda(u)), u \in [0,1].$
Puisque $\{0\} \in \S(t), \forall t \in \C,$ la tressse en question doit \^etre consid\'er\'ee comme {\bf tresse annulaire (annular braid)}, not\'ee $CB_{2m}$
(cylindrical braids), au sens de \cite{KentPeifer}.  C'est une tresse \`a  $2m$ brins correspondant aux $2m$ poin\c{c}ons non-nuls qui font le tour du  segment $\{0\}  \times [0,1] \subset \C \times [0,1].$ Voir figure \ref{fig7}.}
\label{annulaire}
\end{remark}

Puisque chaque racine de $\S(t)$ est une fonction lin\'eaire de la variable complexe $t,$ il est l\'egitime d'\'etendre la notion de rotation \`a celle de rotation fractionnaire. 

\begin{definition} Soit $R_{A_{\J}}$ la rotation d\'efinie dans la d\'efinition \ref{defrotation}.
Pour $\nu = p/q \in [0,1] \cap  \Q, $ on d\'efinit $R_{A_{\J}}^\nu$ comme mouvement des racines  
qui satisfait $ (R_{A_{\J}}^\nu)^q = R_{A_{\J}}^p$ de telle sorte que
$$   Arg\left( \frac{R_{A_{\J}}^\nu (T_{j_{k},\ell_k} (t))  -x }{T_{j_{k},\ell_k} (t)  -x}\right ) =  \nu
Arg \left( \frac{T_{j_{k+1},\ell_{k+1}} (t)  -x }{T_{j_{k},\ell_k} (t)  -x}\right ),   $$
$\exists x \in  ( P_{\J_\mu})^{int} ,$ $\forall k  \in\Z/  |\J |\Z.$
\label{fractionnaire}
Cette d\'efinition peut \^etre \'etendue aux cas de $\nu$ rationnel quelconque si on
se sert de la convention $(R_{A_{\J}})^\nu = (R_{A_{\J}})^{\lfloor\nu \rfloor} \cdot (R_{A_{\J}})^\kappa$
pour $\kappa = \nu -  \lfloor\nu \rfloor.$
\end{definition}

Dans le chapitre pr\'esent, nous ne regardons que des cas avec $\nu = \pm \frac{1}{2},$ tandis que dans  
\cite{TanKocar} nous traitons des cas de rotation \`a puissance rationnelle quelconque des fonctions alg\'ebriques. 

\subsection{Monodromie des tresses induite par des chemins-groupo\"ides}\label{tressesinduites}

Soit la courbe continue $\lambda_{P_i, P_{i+1}}$ sur $\C^\ast \setminus \D$ issue d'un point $P_i$ qui atteint le point $P_{i+1}.$

Il est facile de voir que par la composition
\begin{equation}
\lambda_{P_i, P_{i+2}} = \lambda_{P_i, P_{i+1}}\cdot \lambda_{P_{i+1}, P_{i+2}},
\label{Aii1i2}
\end{equation}
ces courbes  d\'efinissent un groupo\"ide $\pi_1(\C^\ast \setminus \D, B)$ \cite[Example 3, Example 5]{Brown}.  Si l'on
admet l'ensemble de points  base $B$ comme l'espace entier $\C^\ast \setminus \D ,$
alors on obtient  le groupo\"ide fondamental o\`u le chemin joignant deux points est un \'el\'ement de la classe d' homotopie.

%En tenant compte de la Remarque \ref{annulaire}, on peut d\'efinir l'application injective
%du grooupo\"ide de courbes continues sur  $\C^\ast \setminus \D$ \`a celui de 
%rotations fractionnaires (D\'efinition \ref{fractionnaire}).

On appelle $\rho$   l'homomorphisme du groupo\"ide des chemins \`a celui des transformations monodromiques
agissant sur $\S^\ast(\bullet)$.
Pour un chemin $\lambda$ joignant $t_1$ \`a $t_2$ de $\C^\ast \setminus \D,$ on trouve
un diff\'eomorphisme $\rho(\lambda)(\S^\ast(t_1)) = \S^\ast(t_2),$ de telle sorte que l'on ait
\begin{equation}
 \rho:   {\rm chemin}\; \lambda \mapsto  \;  {\rm diff\acute{e}omorphisme}\; \rho(\lambda).
\label{rho}
\end{equation}
Par construction, la famille de diff\'eomorphismes $\rho(\bullet)$ induits par les chemins 
donne naissance \`a un groupo\"ide de transformations monodromiques agissant sur les fibres $\S^\ast(\bullet)$
qui correspond \`a la composition des chemins \eqref{Aii1i2},
$$ \rho(\lambda_{P_i, P_{i+2}}) = \rho(\lambda_{P_i, P_{i+1}})\cdot \rho(\lambda_{P_{i+1}, P_{i+2}}).$$

Une telle famille continue de diff\'eomorphismes induite par le chemin $\lambda$ a \'et\'e utilis\'ee pour d\'ecrire la transformation
de Picard--Lefschetz du cycle \'evanescent provoqu\'ee par un twist de Dehn \cite[1.1]{GZ77}.

L'homomorphisme \eqref{rho} est d\'efini \`a isotopie pr\`es.  Plus pr\'ecisement, si les deux chemins
continus sur $\C^\ast \setminus \D,$
 \begin{equation}
s \in [0,1] \longrightarrow \lambda_\ell = \{t_\ell (s)\} ,\;\; \ell =1,2,
\label{isotope}
\end{equation}
sont homotopes, alors chacune des $2m$ courbes traces $\S^\ast(\lambda_1) = \bigcup_{s \in[0,1]} \S^\ast(t_1 (s)) $ est isotope \`a l'une des courbes traces $\S^\ast(\lambda_2) = \bigcup_{s \in[0,1]} \S^\ast(t_2 (s)). $  Dans ce cas, on dit que $\rho(\lambda_1)$ est isotope \`a  $\rho(\lambda_2).$ 

\begin{remark}
{\rm Il est bien connu que l'on peut indentifier le diff\'eomorphisme $\rho$ agissant sur $\S^\ast (t), t\not \in \D,$ avec un \'el\'ement du groupe de tresses annulaires $CB_{2m}$ (remarque \ref{annulaire}).
L'unique diff\'erence avec le cas classique \cite[corollaire 2.8]{FLP}, concernant le diff\'eomorphisme d'un disque priv\'e de $(2m+1)$ points,
consiste dans l'immobilit\'e de l'origine dans notre cas. }
\label{diffeo}
\end{remark}

Maintenant nous passons \`a la description concr\`ete de la  transformation monodromique \eqref{rho} en fonction du chemin
sous-jacent.

Pour faciliter la formulation des r\'esultats, nous introduisons les r\'eunions des racines en groupes comme suit, qui sont situ\'ees dans $\C^\ast.$

\begin{equation}
A_{j,j}^{(0, 1)} (t) = \{T_{j,0}(t), T_{j,1}(t)\}, \;\; j \in [0; m-1],
\label{Aj01}
\end{equation}

\begin{equation}
 A_{[0 ; m-1]}^{(1)} = \{ T_{j,1}\}_{j=0}^{m-1}, \;\;  resp. \; A_{[0 ; m-1]}^{(0)} = \{ T_{j,0}\}_{j=0}^{m-1}.
\label{A0m-101}
\end{equation}\
Ici les poin\c{c}ons sont rang\'es de telle mani\`ere que 
$ Arg (T_{j, \delta}) < Arg (T_{j+1, \delta}), $ $ (j, \delta) \in [0; m-1] \times [0;1].$
 Pour chaque $\delta \in [0;1]$ fixe, l'enveloppe convexe de $A_{[0 ; m-1]}^{(\delta)}$ contient l'origine. Voir figure \ref{fig7}.

Pour chaque $j \in [0; m-1],$ on regarde le couple comme suit,
\begin{equation}
 A_{j, j+h}^{(0,1)} = \{ T_{j,0}, T_{j+h,1}\},\;\; h \in [0; m-1].
\label{A01jh}
\end{equation}
Ici on utilise la convention $ j+h \in \Z/m\Z.$
En particulier, on a
\begin{equation}
 A_{j, j-k}^{(0,1)} = \{ T_{j,0}, T_{m+j-k,1}\},\;\; k \in [1; m-r-1].
\label{A01jmk}
\end{equation}
L'origine demeure en dehors de l'enveloppe convexe des poin\c{c}ons  $A_{j, j-k}^{(0,1)}$ 
pour chaque $(j,k)$ fixe.

On introduit aussi  la rotation suivante,
 \begin{equation}  
  R_{A_{[+\ell]}^{(0, 1)}} =  \prod_{j=0}^{m-1}R_{A_{j, j+\ell}^{(0, 1)}},
\label{RA+ell01}
\end{equation}
comme la tresse annulaire (remarque \ref{annulaire}) pour $\ell \in [0; r].$
La composition de \eqref{RA+ell01}, $\ell \in [0; k-1],$  produit une concat\'enation des rotations
\begin{equation}  
  R_{A_{[\delta; k-1]}^{(0, 1)}} :=  
  R_{A_{[\delta]}^{(0, 1)}} \cdots  R_{A_{[+(k-1)]}^{(0, 1)}},
\label{RA1k+101}
\end{equation}
pour $\delta \in [0;1].$

D'une mani\`ere similaire, on introduit  
 \begin{equation}  
  R_{A_{[-\ell]}^{(0, 1)}} :=  \prod_{j=0}^{m-1} R_{A_{j, j-\ell}^{(0, 1)}},
\label{RA-ell01}
\end{equation}
pour $\ell \in [1, m-r-1].$

La composition des rotations \eqref{RA-ell01}, $\ell \in [1, k-1],$ entra\^ine la concat\'enation,
\begin{equation}  
  R_{A_{-[1; k-1]}^{(0, 1)}} :=   R_{A_{[-1]}^{(0, 1)}} \cdots R_{A_{[-(k-1)]}^{(0, 1)}}.
\label{RA1k-101}
\end{equation}

Dans notre th\'eor\`eme ci-dessous, qui est l'\'enonc\'e central de ce chapitre,
nous utilisons l'isomorphisme entre le groupe de diff\'eomorphisme $\S^\ast ( - \epsilon \img)$
et le groupe de tresses annulaires $CB_{2m}$ mentionn\'e
dans la remarque \ref{diffeo}.

\begin{theorem}
La liste suivante contient tous les g\'en\'erateurs du groupe de transformations  monodromiques \eqref{rho} pour 
l'arrangement  ${\mathcal A}_{m,2}.$

{\bf (A)} La monodromie autour de l'origine 

\begin{equation} \rho
(\mu_{1/2}^+ \cdot \mu_{2/2}^+ \cdot (\mu_{2/2}^-)^{-1} \cdot (\mu_{1/2}^-)^{-1})=  (  R_{A_{[0]}^{(0, 1)}})^2.
\label{Aj012}
\end{equation}

{\bf (B)} La monodromie autour des points $t = \eta,$ 
\begin{equation}
\rho(\mu^\eta_{1/2}\cdot \lambda_{\eta} \cdot (\mu^\eta_{1/2})^{-1} ) = (  R_{A_{ [0 ; m-1]}^{(\frac{\eta+1}{2})}})^{m},
\label{A0m-1m}
\end{equation}
$\eta = \pm 1.$

{\bf (C)} La monodromie autour du  point $ t= p_{m-k},  k \in [1;m-r-1], $
\begin{equation} (  R_{A_{[-k]}^{(0, 1)}})^2,
\label{A-k012}
\end{equation}
 et celle autour du point $ t= p_k, k \in [1; r],$ \eqref{pk},
\begin{equation} 
 (  R_{A_{[+k]}^{(0, 1)}})^2.
\label{A+k012}
\end{equation}
\label{maintheorem}
\end{theorem}

\begin{proof}
{\bf (A)}
Lorsque $t$ se d\'eplace sur
 l'arc $\mu_{1/2}^-$ joignant $- \epsilon \img$ \`a $  -\epsilon $ (voir figure \ref{fig1}) chacun des $m$ couples 
de racines \eqref{Aj01}
ach\`eve un tour de $- \pi/2$ par rapport au  point $ - \omega^j$ qui peut \^etre represent\'e comme
\begin{equation}
\rho(\mu_{1/2}^- ) = ( R_{A_{[0]}^{(0, 1)}})^{-1/2},
\label{mu12-}
\end{equation}
gr\`ace \`a la commutativit\'e des rotations $ R_{A_{j}^{(0, 1)}},\; j \in [0; m-1].$

D'une mani\`ere similaire, on a,
\begin{equation}
\rho(\mu_{1/2}^- \vec \cdot \mu_{2/2}^-) = ( R_{A_{[0]}^{(0, 1)}})^{-1}.
\label{rhomu12-mu22-}
\end{equation}
La composition des chemins s'effectue  de gauche \`a droite. 
\begin{equation}
\mu_{1/2}^- \vec \cdot \mu_{2/2}^- = \{ - \img e^{ - \img \pi u }; u \in [0, \frac{2}{2}]\}.
\label{mu12-mu22-}
\end{equation}

D'une mani\`ere analogue, on a
\begin{equation}
\rho(\mu_{\ell/2}^+ ) = ( R_{A_{[0]}^{(0, 1)}})^{1/2},
\label{mu12+}
\end{equation}
pour
$$ \mu_{\ell/2}^+ = \{ - \img e^{\img \pi u};   u \in [\frac{\ell -1}{2}, \frac{\ell }{2}  ]\},\;\; \ell = 1,2.  $$
Ici on remarque que le tour complet $\{ - \img e(u);   u \in [0,1 ]\}$ provoque la monodromie de tresses
$$ ( R_{A_{[0]}^{(0, 1)}})^2.$$

\begin{figure}[H]
\centering
\includegraphics[totalheight=10cm]{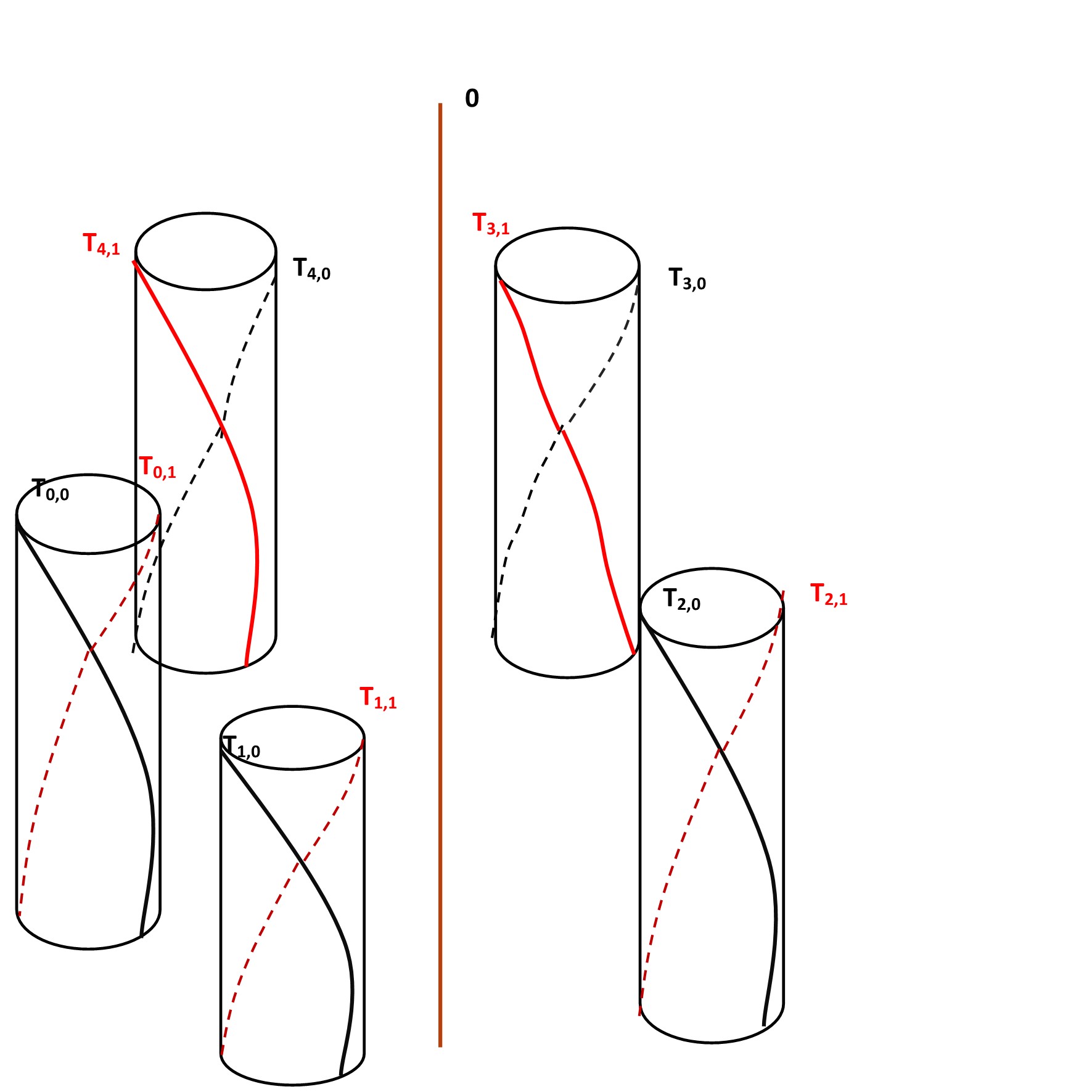}
\caption{ 
Rotation $\rho(\mu_{1/2}^+ \vec \cdot \mu_{2/2}^+) =  \prod_{j=0}^{m-1} R_{A_{j,j}^{(0, 1)}},$ $m=5.$  }
\label{figRAj01}
\end{figure}

{\bf (B)}
Maintenant,  nous examinons la monodromie de tresses induite par le mouvement sur le chemin $\lambda_1:$
d'abord $t$ varie sur $[\epsilon, 1-\epsilon],$ ensuite il fait un tour positif autour $t=1$ le long du cercle
$\{ 1 - \epsilon e(u); u \in [0,1] \}$ et, finalement, il revient \`a $\epsilon$ en suivant  le segment $[1-\epsilon, \epsilon].$

En vertu du lemme \ref{cooincidence} (1), le mouvement le long de $\lambda_1$ (resp. $\lambda_{-1}$) provoque 
la rotation
\begin{equation}
\rho(\lambda_{1}) = (  R_{A_{[0 ; m-1]}^{(1)}})^{m},\;\;  resp. \; \rho(\lambda_{-1}) = (  R_{A_{[0 ; m-1]}^{(0)}})^{m}
\label{lambda1+}
\end{equation}
pour \eqref{A0m-101}.
%$$ A_{[0 ; m-1]}^{(1)} = \{ T_{j,1}\}_{j=0}^{m-1}, \;\;  resp. \; A_{[0 ; m-1]}^{(0)} = \{ T_{j,0}\}_{j=0}^{m-1}.$$

De m\^eme,  le mouvement le long du lacet $\mu^+_{1/2}\cdot \lambda_1 \cdot (\mu^+_{1/2})^{-1}$
provoque la monodromie de tresses
\begin{equation}
\rho(\mu^+_{1/2}\cdot \lambda_1 \cdot (\mu^+_{1/2})^{-1} ) = ( R_{A_{[0]}^{(0, 1)}})^{1/2} (  R_{A_{[0 ; m-1]}^{(1)}})^{m}  ( R_{A_{[0]}^{(0, 1)}})^{-1/2} ,
\label{mu12+lambda1}
\end{equation}
en vertu de \eqref{mu12+}, \eqref{lambda1+}.
D'une mani\`ere similaire, on a
\begin{equation}
\rho(\mu^-_{1/2}\cdot \lambda_{-1} \cdot (\mu^-_{1/2})^{-1} ) = (R_{A_{[0]}^{(0, 1)}})^{- 1/2} (  R_{A_{[0 ; m-1]}^{(0)}})^{m}  ( R_{A_{[0]}^{(0, 1)}})^{1/2}.
\label{mu12-lambda-1}
\end{equation}

On voit facilement que les rotations  $(R_{A_{[0]}^{(0, 1)}})^{\pm \frac{1}{2}},$ et  $ R_{A_{[0 ; m-1]}^{(0)}}$ commutent.
Pour le voir, nous consid\'erons le tuyau $B_j$ de diam\`etre assez petit qui contient $A_{j,j}^{(0, 1)}$ \eqref{Aj01} et nous notons
\begin{equation}
B_{[0;m-1]} = \cup_{j=0}^{m-1} B_j.
\label{B0m-1}
\end{equation} 
La rotation $R_{B_{[0;m-1]}}$ est isotope \`a la rotation $ R_{A_{[0 ; m-1]}^{(0)}}.$
Puisque la rotation $R_{A_{j,j}^{(0, 1)}}$ a lieu \`a l'int\'erieur du tuyau  $B_j,$ elle commute avec $  R_{A_{[0 ; m-1]}^{(0)}}.$

\begin{figure}[H]
\centering
\includegraphics[totalheight=13cm]{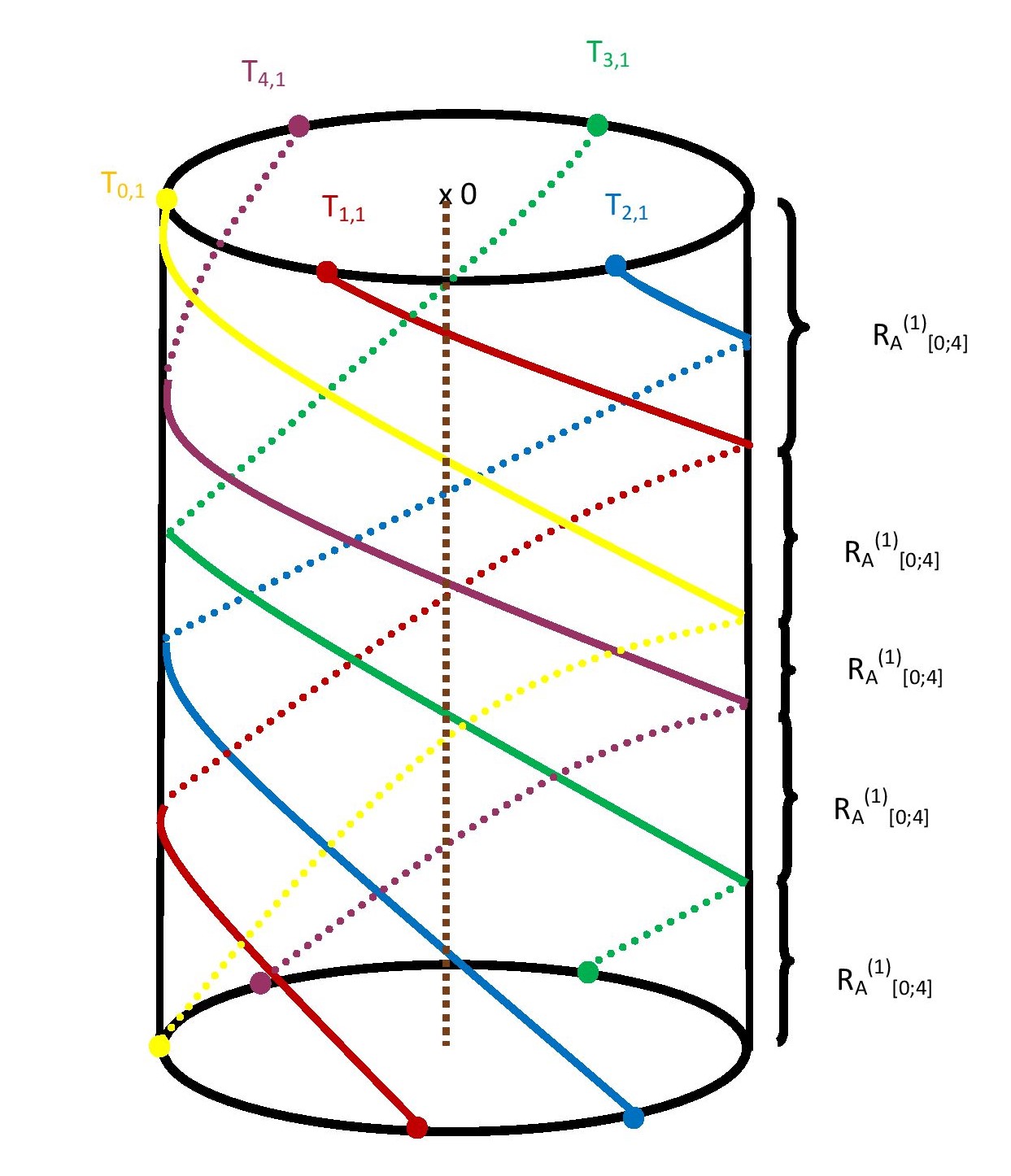}
\caption{ Rotation $(  R_{A_{[0 ; 4]}^{(1)}})^{5}.$  }
\label{fig7}
\end{figure}

{\bf (C) }
Maintenant, nous examinons la monodromie de tresses induite par le mouvement de $t = -\epsilon \img$ au point de co\"incidence,
\begin{equation}
p_{m-k} = \img \tan (\frac{m-k}{m}) = -  \img \tan (\frac{k}{m}),\;\; k \in [1; r+1]
\label{pm-k}
\end{equation}
avec
$$r := \lfloor \frac{m-1}{2}\rfloor.$$

 Le mouvement qui part de  $t = - \epsilon \img $ et arrive \`a $ t =  \epsilon \img + p_{m-1},$ en suivant  au d\'ebut une droite, ensuite un tour sur le cercle $ p_{m-1} + \epsilon e(u), u \in [0;1],$ induit la rotation
 \begin{equation}  
 ( R_{A_{[-1]}^{(0, 1)}})^2= \prod_{j=0}^{m-1}( R_{A_{j, j-1}^{(0, 1)}})^2,
\label{Rotpm1}
\end{equation}
puisqu'au point $ t =  p_{m-1},$ les deux poin\c{c}ons de $ A_{j, j-1}^{(1)}$ \eqref{A01jmk} se rencontrent.

\begin{figure}[H]
\centering
\includegraphics[totalheight=10cm]{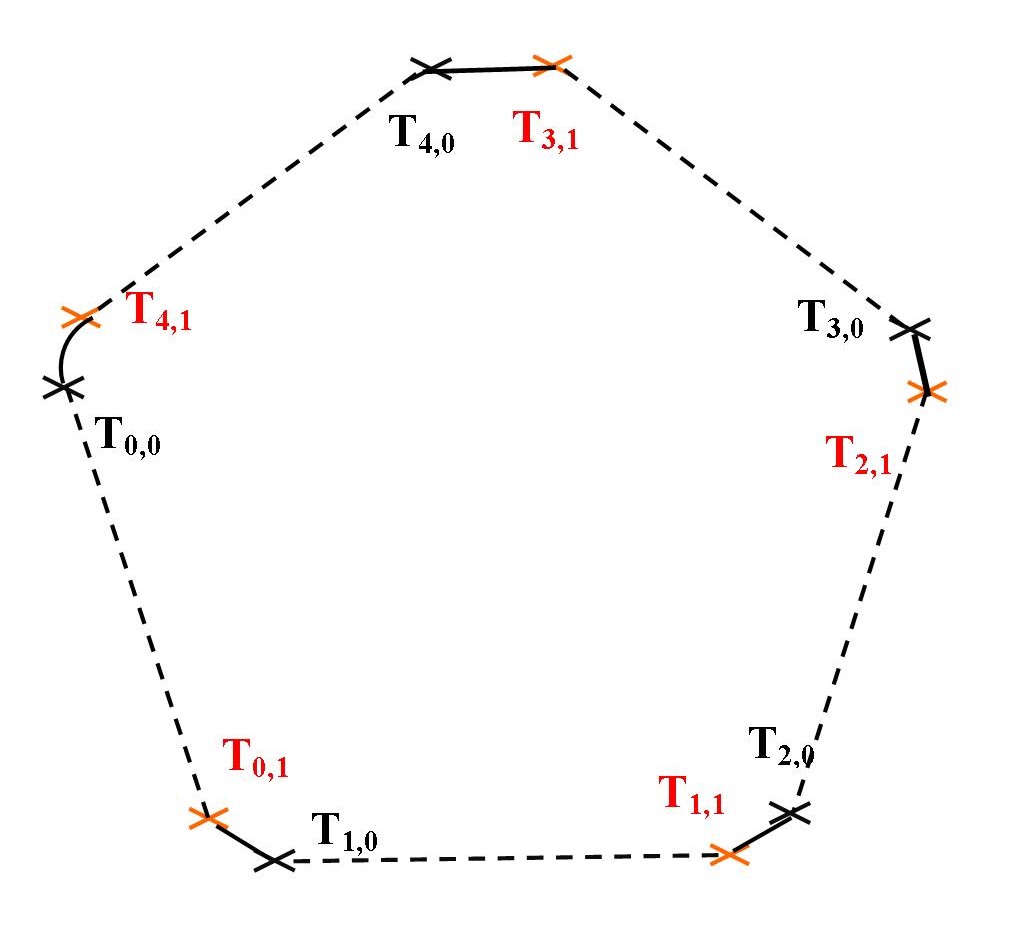}
\caption{ $\S( \epsilon \img + p_{m-1}), m=5.$}
\label{fig8}
\end{figure}

De la figure \ref{fig8}, on voit que les rotations $R_{A_{j, j-1}^{(0, 1)}},  j \in [1; m-1],$ commutent entre elles.
%D'une mani\`ere similaire, on voit la commutativit\'e des rotations $R_{A_{j, j-\ell}^{(0, 1)}},  j \in [1; m-1],$
%pour $\ell$ fixe. 

Pour examiner la monodromie provoqu\'ee en $t = p_{m-k}, k \in [2;r+1], $ il faut consid\'erer les demi-tours
au sens, d'abord  positif et au retour ensuite n\'egatif,  pr\`es des points   $t = p_{m-\ell}, \ell \in [1; k-1]$ comme on le voit
dans la figure \ref{fig1}, o\`u l'on a omis le mouvement de retour.\\

En  $t = p_{m-\ell}, \ell \in [1; k-1],$
lors du demi-tour en sens positif $ p_{m-\ell} + \epsilon e(u), u \in [0;1/2],$ 
les racines $\S^\ast(t)$ \eqref{sectionSt} subissent la rotation 
$  R_{A_{[-\ell]}^{(0, 1)}} $ \eqref{RA-ell01}.

\begin{figure}[H]
\centering
\includegraphics[totalheight=25cm]{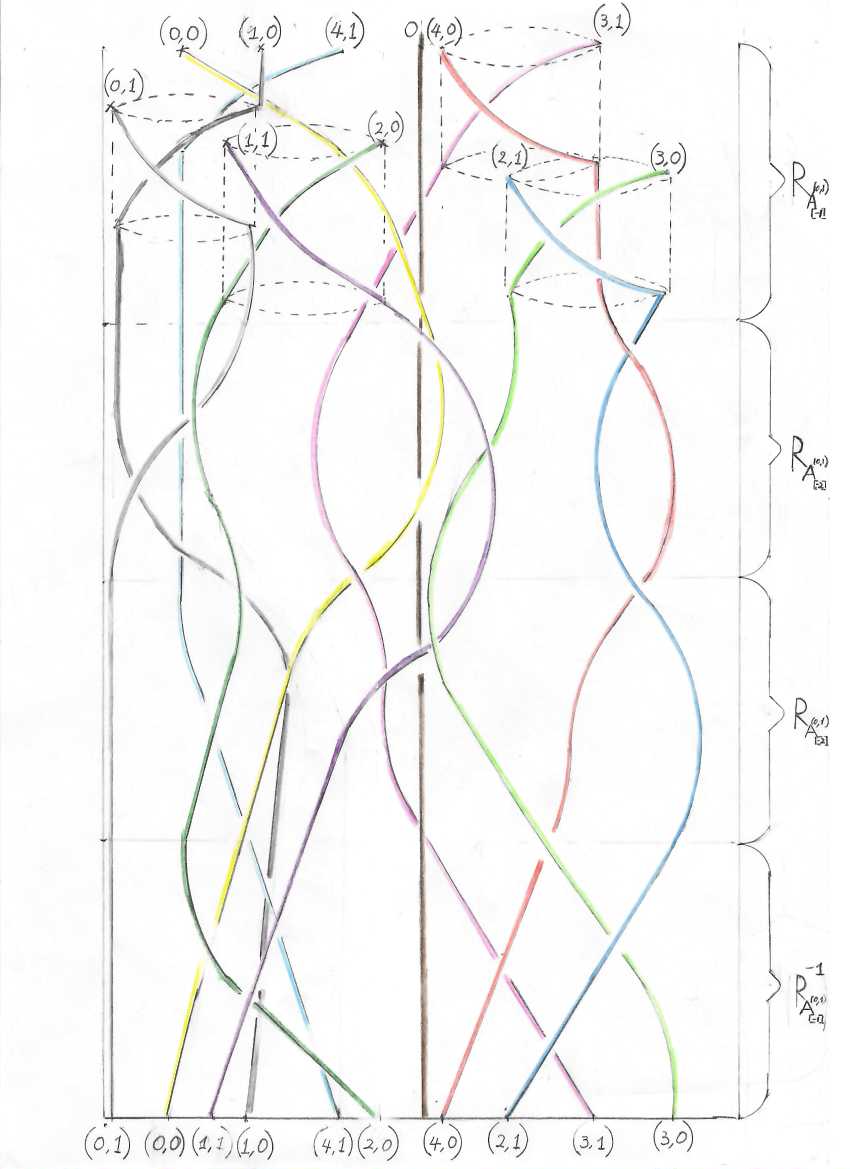}
\caption{ 
$R_{A_{[-1]}^{(0, 1)}} \cdot (  R_{A_{[-2]}^{(0, 1)}})^2  \cdot R_{A_{[-1]}^{(0, 1)}}^{-1}$ pour $m=5.$ }
\label{fig8bis}
\end{figure}

Le d\'etour positif de tous les points $t = p_{m-\ell}, \ell \in [1; k-1]$  donne naissance \`a
la composition
$  R_{A_{-[1; k-1]}^{(0, 1)}}$
\eqref{RA1k-101}.

Le tour complet autour de $ t= p_{m-k}$ \eqref{RA1k-101}
produit  la rotation
 \begin{equation}  
(  R_{A_{[-k]}^{(0, 1)}})^2 =  \prod_{j=0}^{m-1} (R_{A_{j, j-k}^{(0, 1)}})^2,
\label{RA-k01}
\end{equation}
 en vertu de la d\'efinition \eqref{RA-ell01} appliqu\'ee au cas $\ell =k.$

Le d\'etour n\'egatif de tous les points $t = p_{m-\ell}, \ell \in [1; k-1],$ induit  $R_{A_{-[1; k-1]}^{(0, 1)}}^{-1}.$

Ainsi le lacet du point base $t = - \epsilon \img$ qui fait un tour positif autour de  $ t= p_{m-k},$ suivant le chemin de la figure \ref{fig1}, provoque
la monodromie de tresses
 \begin{equation}  
R_{A_{-[1; k-1]}^{(0, 1)}} \cdot (  R_{A_{[-k]}^{(0, 1)}})^2  \cdot R_{A_{-[1; k-1]}^{(0, 1)}}^{-1}.
\label{Rotpk01}
\end{equation}

La figure \ref{fig8bis} montre que la rotation \eqref{Rotpk01} est isotope \`a la rotation  $(  R_{A_{[-k]}^{(0, 1)}})^2$
pour le cas $k=2.$

 De mani\`ere inductive, on montre  la r\'eduction de \eqref{Rotpk01} \`a   \eqref{A-k012} pour $k$ g\'en\'eral. 
L'argument s'appuie sur l'observation suivante. Lors de la rotation  \eqref{Rotpk01}, les brins $T_{j-h, \delta}, (h, \delta) \in [1; k-1] \times [0;1],$
restent devant ou derri\`ere les brins $A^{(0,1)}_{j, j-k}$ \eqref{A01jmk} pour chaque $j \in [0; m-1].$
C'est-\`a-dire,
dans \eqref{Rotpk01} la tresse $R_{A^{(0,1)}_{j, j-k}}$ n'est pas conjugu\'ee avec des tresses affect\'ees par $T_{j-h, \delta}, (h, \delta) \in [1; k-1] \times [0;1].$
Les autres brins $T_{j-\ell, \delta},$ $ (\ell, \delta) \in [k+1; m-1]\times [0;1],$ se trouvent aussi \'eloign\'es des brins
$T_{j-h, \delta}, (h, \delta) \in [0; k] \times [0;1],$ qu'elles ne s'entrelacent avec les derni\`eres.
En se servant de la notation r\'epandue de \cite{KentPeifer}, nous pourrions r\'esumer cette relation par la formule
$$ \sigma_i (\sigma_j)^2 (\sigma_i)^{-1} =  (\sigma_j)^2, \;\; |i - j| \not = 1, 2m-1,$$
pour $\sigma_i, \sigma_j \in CB_{2m}.$

En s'appuyant sur un argument similaire, on peut aussi r\'eduire \eqref{mu12-lambda-1} \`a la rotation $(  R_{A_{[0 ; m-1]}^{(0)}})^{m}.$

Quant \`a la monodromie de tresses provoqu\'ee par le tour complet autour du point $t = p_{k}, k \in [1; r],$
sa description est parall\`ele \`a celle pour les cas $ k \in [r+1; m-1]$ donn\'ee dans \eqref{Rotpk01}.

La rotation $\rho(\mu_{1/2}^+\vec \cdot \mu_{2/2}^+)$ \eqref{mu12+}
provoqu\'ee par  le demi-tour autour de l'origine
est  suivie des rotations induites par les demi-tours  $R_{A_{[+\ell]}^{(0, 1)}}$ \eqref{RA+ell01}  autour les points $t = p_{\ell}, \ell \in [1; k-1],$
dont la composition s'\'ecrit comme 
$R_{A_{[1; k-1]}^{(0, 1)}}$ 
\eqref{RA1k+101}.

Le tour complet  autour de $ t= p_k,$ \eqref{RA+ell01}
produit  la rotation
 \begin{equation}  
(  R_{A_{[+k]}^{(0, 1)}})^2 =  \prod_{j=0}^{m-1} (R_{A_{j, j+k}^{(0, 1)}})^2.
\label{RA+k012}
\end{equation}

Ainsi le lacet du point base $t = - \epsilon \img$ qui fait un tour positif autour de  $ t= p_k, k \in [1; r]$ suivant le chemin de la figure \ref{fig1}, provoque
la monodromie de tresses
 \begin{equation}  
 R_{A_{[0; k-1]}^{(0, 1)}} \cdot (  R_{A_{[+k]}^{(0, 1)}})^2  \cdot R_{A_{[0; k-1]}^{(0, 1)}}^{-1} .
\label{Rotp+k01}
\end{equation}
pour $k \in [1;r].$

 Pour voir que la rotation
 \begin{equation}  
 (R_{A_{[0]}^{(0, 1)}})\cdot (  R_{A_{[+k]}^{(0, 1)}})^2  \cdot (R_{A_{[0]}^{(0, 1)}})^{-1} .
\label{Rotp+k01simple}
\end{equation}
est isotope \`a \eqref{A+k012}, on a recours \`a l'examen qui est similaire
\`a la r\'eduction de \eqref{Rotpk01} \`a   \eqref{A-k012}.

Les figures \ref{figRA1SqConj}, \ref{figRA1Sq} illustrent le fait que la rotation
$(R_{A_{[0]}^{(0, 1)}})\cdot (  R_{A_{[+k]}^{(0, 1)}})^2  \cdot (R_{A_{[0]}^{(0, 1)}})^{-1}$ est isotope \`a $(  R_{A_{[+k]}^{(0, 1)}})^2 $ pour $k=1.$

\begin{figure}[H]
\centering
\includegraphics[totalheight=20cm]{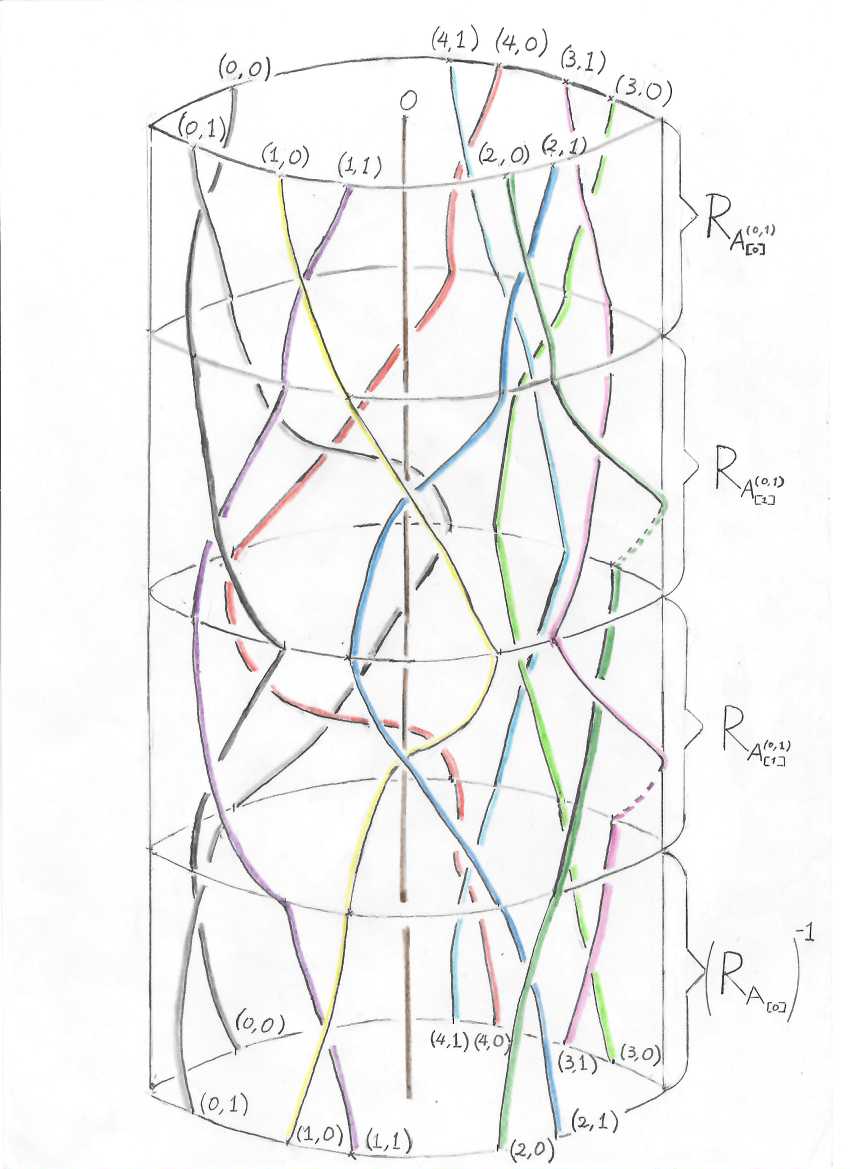}
\caption{ Rotation $(R_{A_{[0]}^{(0, 1)}})\cdot (  R_{A_{[1]}^{(0, 1)}})^2  \cdot (R_{A_{[0]}^{(0, 1)}})^{-1}$ pour $m=5.$}
\label{figRA1SqConj}
\end{figure}

\begin{figure}[H]
\centering
\includegraphics[totalheight=18cm]{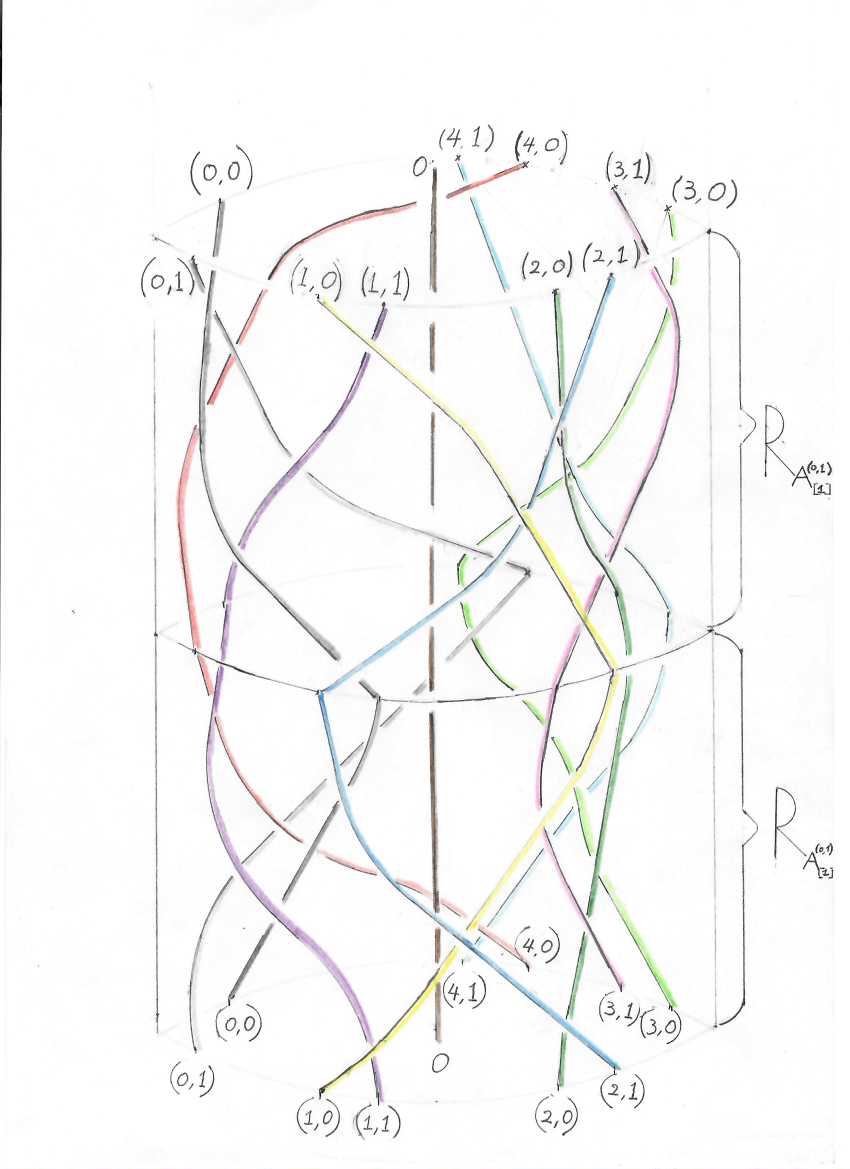}
\caption{ Rotation $(  R_{A_{[1]}^{(0, 1)}})^2 $ pour $m=5.$}
\label{figRA1Sq}
\end{figure}

%La sym\'etrie entre deux rotations, \eqref{Rotpk01} et \eqref{Rotp+k01}, est \'evidente.

\end{proof}

Pour \'etablir une description du groupe fondamental 
$\pi_1(\C^2 \setminus  {\mathcal A}_{m,2}),$ nous consid\'erons les lacets
\begin{equation}  
T= \langle \tau_{j,0}, \tau_{j,1} \rangle_{j=0}^{m-1} \cup\{\tau_0\}
\label{tau}
\end{equation}
qui forment un syst\`eme des chemins distingu\'es \cite[Definition 2]{Gabrielov73}, \cite[Definition 1.2.3]{GZ77},\cite[Definition 1.2]{Loenne} de $\S(- \epsilon \img).$ Le lacet $\tau_{j,\delta},$
partant d'un point fixe, 
fait un tour positif autour du poin\c{c}on $T_{j,\delta}(- \epsilon \img),$ $(j, \delta) \in [0; m-1] \times [0;1].$

\begin{figure}[H]
\centering
\includegraphics[totalheight=8cm]{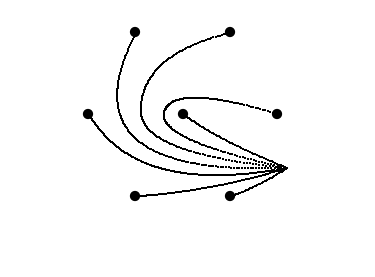}

\caption{  Chemins distingu\'es pour $m=3.$}
\label{distinguished}
\end{figure}

Consid\'erons maintenant l'action de la transformation monodromique $\rho$ sur les lacets $\tau_{j,\delta}$
$\rho_\ast(\tau_{j,\delta}).$

\begin{figure}[H]
\centering
\includegraphics[totalheight=8cm]{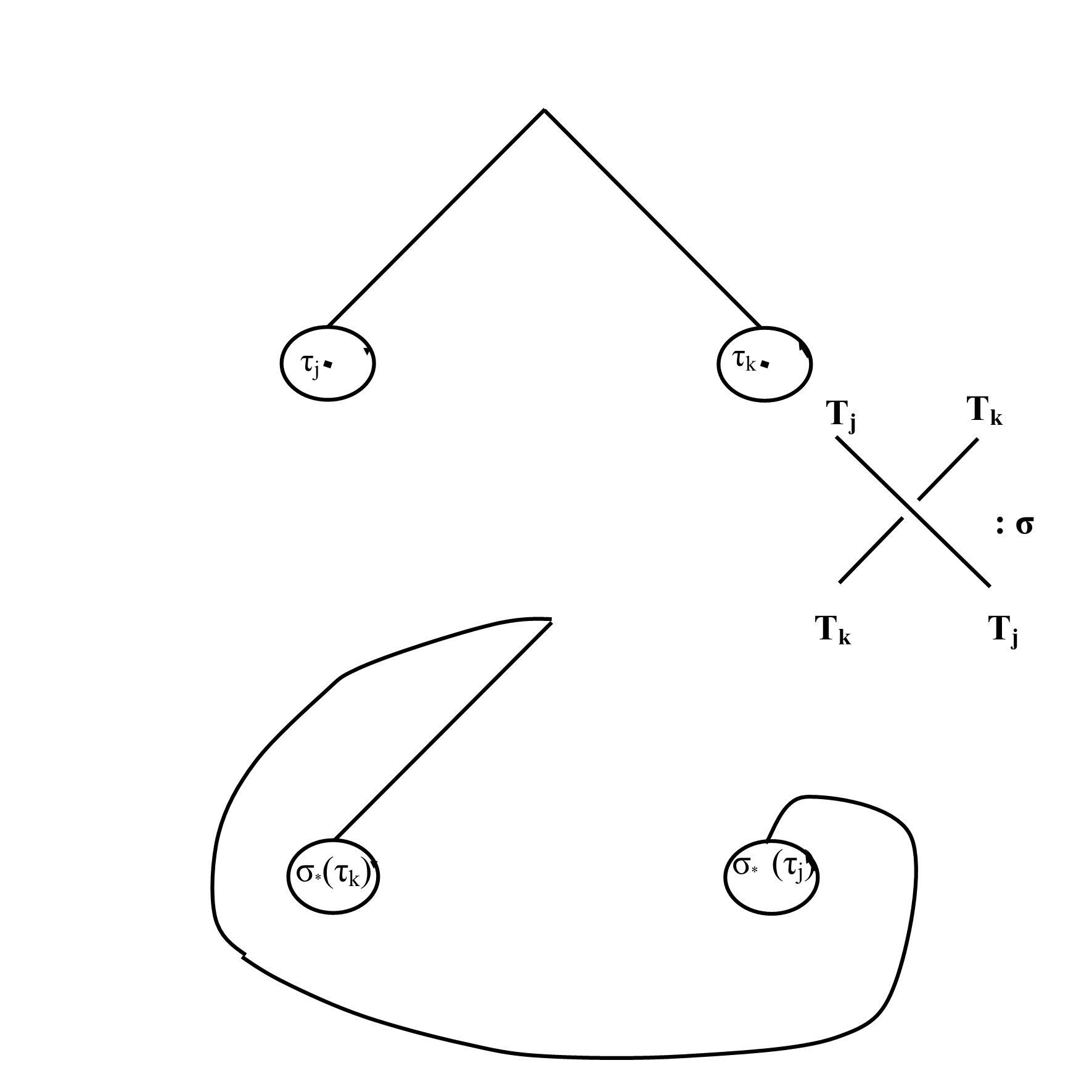}

\caption{ Action de la tresse $\sigma$ sur les lacets $\tau_j, \tau_\ell$
: $\sigma_\ast(\tau_j) = \tau_j \tau_\ell {\tau_j}^{-1}.$ }
\label{braidloop}
\end{figure}

Nous pouvons, alors, d\'eduire le corollaire suivant du th\'eor\`eme 
\ref{maintheorem}.
\begin{corollary}
On a la description suivante du groupe fondamental du compl\'ement de l'arrangement ${\mathcal A}_{m,2}$ en termes
des g\'en\'erateurs $T$ \eqref{tau},
\begin{equation} 
\pi_1(\C^2 \setminus  {\mathcal A}_{m,2}) = \{ T; \tau_{j,\delta} = \rho_\ast(\tau_{j,\delta} ), (j, \delta) \in [0; m-1] \times [0;1] \}.
\label{pi1Am2}
\end{equation}
Ici $\rho$ est le diff\'eomorphisme de $\S^\ast( - \epsilon \img)$ induit par une des rotations \eqref{Aj012}, \eqref{A0m-1m},  \eqref{A-k012}, \eqref{A+k012} de la liste du
th\'eor\`eme  \ref{maintheorem}. 
%{\color{blue} On remarque que $\tau_0$ et  $\rho_\ast(\tau_{j,\delta} )$ \eqref{pi1Am2} commutent pour $\rho$ induit par la rotation \eqref {Aj012}, mais ne commutent pas pour $\rho$ induit par \eqref{A0m-1m}, \eqref{A-k012} et \eqref{A+k012}.}
\label{fgroup}
\end{corollary}
 
 La d\'emonstration n'est qu'une application directe du lemma 1.14 de \cite{Loenne} (th\'eor\`eme de Zariski--Van Kampen) \`a notre situation.  

 \vspace{1pc}

En conclusion, nous proposons la conjecture suivante vu la r\'egularit\'e que l'on observe dans l'expression des rotations 
\eqref{Aj012}, \eqref{A0m-1m}, \eqref{A-k012}, \eqref{A+k012}.
Notamment, le th\'eor\`eme \ref{maintheorem} implique que   la monodromie globale pour l'arrangement ${\mathcal A}_{m,2}$ est compl\`etement d\'etermin\'ee par la monodromie locale. 

Pour formuler la conjecture, nous introduisons quelques notations plus g\'en\'erales que celles du cas $n=2.$

Soient les racines $s(t)$ de \eqref{Amn},
\begin{equation}
T_{\bk_\alpha} (t)=- \omega_m^{j_\alpha} ( 1+\omega_n^{\ell_\alpha} t )
\label{Tjlmn}
\end{equation}
pour $\bk_\alpha = ({j_\alpha}, \ell_\alpha) \in [0; m-1] 
\times [0;n-1].$
On note l'ensemble
des points de co\"incidence  $\D_{m,n}$ avec
\begin{equation}
   \D_{m,n} =\{ t \in \C; T_{\bk_\alpha} (t)= T_{\bk_\beta} (t), \exists \bk_\alpha \not = \bk_\beta \in [0; m-1] 
\times [0;n-1] \}.
\label{coincidenceDmn}
\end{equation} 
Pour $\xi \in \D_{m,n}$ fixe, on d\'efinit les ensembles d'indices $\{ \Lambda_{\mu} (\xi)\}_{\mu=1}^{\nu_\xi}$
tels que
\begin{equation}
   T_{\bk_\alpha} (\xi)= T_{\bk_\beta} (\xi), \forall \bk_\alpha \not = \bk_\beta \in \Lambda_{\mu} (\xi) \subset [0; m-1] \times [0;n-1] 
\label{coincidenceDmnT}
\end{equation} 
pour $\mu \in [1; \nu_\xi].$
Ici on suppose que $ T_{\bk} (\xi)\not = T_{\bh} (\xi), $ $\forall \bk \in \Lambda_{\mu} (\xi),$ $\forall \bh \not \in \Lambda_{\mu} (\xi)$ pour chaque $\mu \in [1; \nu_\xi].$
Dans l'ensemble $\Lambda_{\mu} (\xi),$ les poin\c{c}ons sont rang\'es selon l'ordre de croissance de l'argument comme dans la d\'efinition \ref{defrotation}.

Le r\'esultat principal de \cite{TanLonne} sera bas\'e sur l'\'enonc\'e de la conjecture \ref{conj} ou sur une description analogue du groupe de monodromie globale.
  
\begin{conjecture}
{\rm Le groupe de transformations  monodromiques globales
pour l'arrangement  ${\mathcal A}_{m,n},$ $n \geq2,$
est compl\`etement d\'etermin\'e par ses monodromies locales des fonctions lin\'eaires \eqref{Tjlmn}. 

Autrement dit, lorsque le param\`etre $t$ 
varie le long d'un lacet issu d'un point $ p \in \C \setminus \D_{m,n}$ qui fait le tour autour de $\xi$ et pas d'autres points de $\D_{m,n}$   
\eqref{coincidenceDmn},
les poin\c{c}ons  $\{T_{\bk} (t)\}_{\bk \in \Lambda_{\mu} (\xi)},$ $\mu \in [1; \nu_\xi],$ \eqref{coincidenceDmnT}
subissent la transformation monodromique
$$ (R_{A_\mu(\xi)})^{|\Lambda_{\mu} (\xi)|}.$$}
\label{conj}
\end{conjecture}
En particulier, on remarque que
$$ |\Lambda_{1} (0)| = \cdots = |\Lambda_{m} (0)| =n$$
$\forall (m,n) \in (\Z_{\geq 2})^2.$

\vspace{\fill}

\begin{flushleft}
 \begin{minipage}[t]{8cm}
  \begin{center}
{\footnotesize 
%Kumamoto University\\
%Kurokami 2-39-1,\\
 %Kumamoto,860-8555,
%Japan\\

%Galatasaray University\\
%\c{C}{\i}ra$\rm\breve{g}$an cad. 36,\\
%Be\c{s}ikta\c{s}, Istanbul, 34357,
%Turkey\\

Moscow Institute of Physics and Technology (National Research University)\\
 Department of Discrete Mathematics, \\
Dolgoprudny, Moscow Region, Russian Federation

{\it E-mail}: tanabesusumu@hotmail.com, tanabe.s@mipt.ru}
\end{center}
\end{minipage}\hfill
\end{flushleft}

\end{document}